\documentclass[a4paper,twoside,10pt]{amsart}

\usepackage[utf8x]{inputenc}

\usepackage{amssymb}
\usepackage{amsmath}
\usepackage{epsfig}
\usepackage{todonotes}

\usepackage[T1]{fontenc}

\newtheorem{The}{Theorem}[section]
\newtheorem{Lem}[The]{Lemma}

\newtheorem{Pro}[The]{Proposition}

\newtheorem{Rem}[The]{Remark}
\newtheorem{Defi}[The]{Definition}

\newtheorem{Exam}[The]{Example}
\newtheorem{Not}[The]{Notation}

\def\C{\mathbb C}

\def\l{\lambda}

\def\N{\mathbb N}

\def\Q{\mathbb Q}
\def\R{\mathbb R}

\def\Z{\mathbb Z}

\date{\today}

\title{On certain polynomial systems involving Stirling numbers of second kind}

\author{F. J. Castro-Jim\'enez}

\address{Departamento de Álgebra e IMUS, Universidad de Sevilla}

\email{castro@us.es}

\author{H. Cobo Pablos}

\email{helenacobo@gmail.com}

\thanks{Both authors are partially supported by MTM2016-75024-P and Feder, and FQM333}
\begin{document}

\begin{abstract}
We solve a special type of linear systems with coefficients in multivariate polynomial rings. These systems arise in the computation of parametric Bernstein-Sato polynomials associated with certain hypergeometric ideals in the Weyl algebra.
\end{abstract}

\subjclass[2010]{15A54; 11B73}

\keywords{Linear systems with polynomial coefficients. Stirling numbers of second kind}

\maketitle

\section{Introduction}
\label{StInt}

The Stirling numbers of first kind, denoted $s(d,k)$, for $d\geq k$, can be defined by the equation
\[x(x-1)(x-2)\cdots(x-d+1)=\sum_{k=0}^ds(d,k)x^k.\]
These well-known combinatorial numbers could also be defined, in a very fancy way, as the solution to the system
\vspace{2mm}
\begin{equation}
\left(\begin{array}{cccc}
S(0,0) & S(1,0) & \cdots & S(d-1,0)\\
0 & S(1,1) & \cdots & S(d-1,1)\\
  &  & \ddots & \\
0 & & 0 & S(d-1,d-1)\\
\end{array}\right)\left(\begin{array}{c}
a_0\\
a_1\\
\vdots\\
a_{d-1}\\
\end{array}\right)=-\left(\begin{array}{c}
S(d,0)\\
S(d,1)\\
\vdots\\
S(d,d-1)\\
\end{array}\right)
\label{SystemS1}
\end{equation}
where $S(i,j)$ denotes the so called Stirling numbers of second kind, another combinatorial numbers closely related to the first kind ones. This is due to the existent relations between the Stirling numbers of first and second kind.
Notice that there is an obvious correspondence between the rows of the $d\times d$ matrix $\big(S(i,j)\big)_{i,j}$ and the points $\{0,1,\ldots,d-1\}$, which are the roots of the polynomial $x^d+a_{d-1}x^{d-1}+\cdots+a_0$.

\vspace{2mm}

In this paper we study a special kind of linear systems which may be seen as a generalization of the system in (\ref{SystemS1}). More precisely, we define a special kind of linear systems whose coefficients are polynomials involving Stirling numbers of second kind.
\vspace{2mm}

Given a set of $r$ points $R\subseteq\Z^2$, we will define, for every point in $R$, a linear equation in the variables $a_0,\ldots,a_{r-1}$ with polynomial coefficients in $\C[x,y,z]$. In this way we associate to $R$ a linear system of equations. Can we solve these systems? The answer is yes under certain conditions on the set $R$. It turns out that, considering the polynomial $b_{\bf a}(s):=s^r+a_{r-1}s^{r-1}+\cdots+a_1s+a_0$ encoding the solution ${\bf a}=(a_0,\ldots,a_{r-1})$ to the system, the roots of $b_{\bf a}(s)$ are also closely related with the points in $R$. We also study this kind of systems in rings of the form $\C[x,y,z]/(ax+by)$ with $(a,b)\in\Z^2$.

\vspace{2mm}

This type of systems (more precisely, specializations of this type of systems to systems with complex coefficients) appear when computing the global $b$-function of hypergeometric systems $H_A(\beta)$ with matrix $A\in\mathcal M_{n,n+2}$. The Stirling numbers of second kind appear naturally in this context since the hypergeometric systems are defined in the Weyl algebra (see \cite{MS}), where the following relation holds
\[\big(x\partial_x)^n=\sum_{k=0}^nS(n,k)x^k\partial_x^k.\]
In \cite{CC} we use the results of this paper to compute de $b$--function of $H_A(\beta)$ where $A$ belongs to a family of matrices of size $1\times 3$ and $\beta$ belogs to a certain set of complex parameters. Some results of this paper can be also applied in the computation of the $b$--function associated with hypergeometric systems for matrices of size $n\times (n+2)$, $\beta \in \C^{n}$ and $n\geq 2$.

\vspace{2mm}

As a consequence of the combinatorics we have developed in the main sections of the paper, we deduce in Section \ref{Appl} a family of relations on the Stirling numbers of second kind. We do not know any reference for this family of relations.

\vspace{3mm}

{\bf Acknowledgments:} We wish to thank Christian Krattenthaler for providing the proof of the essential Proposition \ref{determinant}.

\section{On Stirling numbers}

A monic polynomial in $\C[s]$ is completely determined either by its coefficients or by its roots. The relation between roots and coefficients is easily understood by using elementary symmetric functions. Let us set some notation.

\vspace{3mm}

\begin{Not}
 Given a vector ${\bf a}=\big(a_0,\ldots,a_{r-1}\big)\in\C^r$ we define the monic polynomial of degree $r$
\[b_{{\bf a}}(s)=s^r+a_{r-1}s^{r-1}+\cdots+a_0.\]
Given a set of $r$ complex numbers $\Delta=\{\alpha_1,\ldots,\alpha_r\}$, we define the monic polynomial of degree $r$
\[b_\Delta(s)=\prod_{i=1}^r\big(s-\alpha_i\big).\]
The relation between ${\bf a}$ and $\Delta$ is determined by the elementary symmetric functions. That is, for $0\leq j\leq r-1$, the coefficients of the polynomial $b_{\Delta}(s)$ are given by
\[a_j=\sum_{1\leq i_1<i_2<\ldots<i_{r-j}\leq r} \alpha_{i_1}\cdots \alpha_{i_{r-j}},\]
the $j$-th elementary symmetric function on the elements of $\Delta$.
\label{Not2}
\end{Not}
\vspace{2mm}

Let us look at a very simple example. Consider the set $\Delta=\{0,1,\ldots,d-1\}$. Then the coefficients of the polynomial $b_\Delta(s)$ are the well known Stirling numbers of first kind.

These combinatorial numbers, as well as the Stirling numbers of second kind, appear in very different branches of mathematics. There are many different ways to define them, as we see next (we refer to \cite{Kr},\cite{C} and \cite{GKP}, which contain very good introductions to Stirling numbers of first and second kind).

\vspace{2mm}

$\bullet$ Stirling numbers of the first kind, denoted $s(d,k)$, can be defined by the generating function:
\begin{equation}
x(x-1)(x-2)\cdots(x-d+1)=\sum_{k=0}^ds(d,k)x^k
\label{eq1S}
\end{equation}
with the conventions $s(n,0)=0$ for $n>0$, $s(n,m)=0$ if $n<m$ and $s(0,0)=1$.

\vspace{2mm}

$\bullet$ Stirling numbers of the second kind, denoted by $S(d,k)$, can be defined by the generating function:
\begin{equation}
x^d=\sum_{k=0}^dS(d,k)x(x-1)(x-2)\cdots(x-k+1)
\label{eqS2}
\end{equation}
with the conventions $S(n,0)=0$ for $n>0$, $S(n,m)=0$ if $n<m$ and $S(0,0)=1$.

\vspace{3mm}

The Stirling numbers of first and second kind are somewhat inverse to each other, since we have
\vspace{2mm}
\begin{equation}
\begin{array}{c}
\left(\begin{array}{cccc}
S(0,0) & S(1,0) & \cdots & S(d,0)\\
0 & S(1,1) & \cdots & S(d,1)\\
  &  & \ddots & \\
0 & & 0 & S(d,d)\\
\end{array}\right)
\left(\begin{array}{cccc}
s(0,0) & s(1,0) & \cdots & s(d,0)\\
0 & s(1,1) & \cdots & s(d,1)\\
  &  & \ddots & \\
0 & & 0 & s(d,d)\\
\end{array}\right)=\\
\\
=\left(\begin{array}{cccc}
s(0,0) & s(1,0) & \cdots & s(d,0)\\
0 & s(1,1) & \cdots & s(d,1)\\
  &  & \ddots & \\
0 & & 0 & s(d,d)\\
\end{array}\right)\left(\begin{array}{cccc}
S(0,0) & S(1,0) & \cdots & S(d,0)\\
0 & S(1,1) & \cdots & S(d,1)\\
  &  & \ddots & \\
0 & & 0 & S(d,d)\\
\end{array}\right)
=I_{d+1}\\
\end{array}
\label{Ssinv}
\end{equation}
%\vspace{2mm}
where $I_{d+1}$ denotes the identity matrix of order $d+1$. Equation (\ref{Ssinv}) can be written in compact form as the following set of relations (orthogonality property):
\begin{equation}
\sum_{k=m}^ds(k,m)S(d,k)=\sum_{k=m}^ds(d,k)S(k,m)=\delta_{md}
\label{InvStr}
\end{equation}
where $\delta_{md}$ denotes the Kronecker delta. This is the clue to solve the system (\ref{SystemS1}).

\vspace{2mm}

Stirling numbers of first and second kind satisfy the following recurrences:
\[s(n+1,k)=s(n,k-1)+ns(n,k)\]
\[S(n+1,k)=S(n,k-1)+S(n,k)\]
respectively. Moreover there are closed forms to compute them:
\begin{equation}
s(n,m)=\sum_{k=0}^{n-m}(-1)^k\binom{n-1+k}{n-m+k}\binom{2n-m}{n-m-k}S(n-m+k,k),
\label{S1}
\end{equation}
\begin{equation}
S(n,m)=\frac{1}{m!}\sum_{k=0}^m(-1)^{m-k}\binom{m}{k}k^n.
\label{S2}
\end{equation}
Both numbers have a combinatorial interpretation. The number $(-1)^{n-m}s(n,m)$ is the number of permutations of $n$ symbols which have exactly $m$ cycles, while the Stirling number of second kind $S(m,n)$ counts the number of ways to partition a set of $m$ elements into $n$ nonempty subsets.

\vspace{3mm}

We end this section with some interesting relations concerning Stirling numbers of second kind.
First a very useful generating function, which can be indeed taken as another definition for Stirling numbers of second kind,
\begin{equation}
\frac{1}{n!}\big(e^x-1\big)^n=\sum_{m=n}^\infty S(m,n)\frac{x^m}{m!}.
\label{FormulaKT}
\end{equation}
Moreover, we have the following well-known identity (see \cite{GKP})
\[k^n=\sum_{m=1}^k\binom{k}{m}m!S(n,m).\]
Finally we state a nice formula (see \cite{Kr}):
\begin{equation}
\sum_{k=0}^n\binom{n}{k}(-1)^k\big(xk+y)^m=(-1)^nn!\sum_{j=n}^m\binom{m}{j}x^jy^{m-j}S(j,n).
\label{eqPalma}
\end{equation}
In Section \ref{Appl} we will give a generalization of this equation.

\vspace{3mm}

\section{On certain polynomial systems of equations in $\C[x,y,z]$}
\label{S2}
In this section we define a particular type of linear systems with polynomial coefficients, and study its solutions in the ring $\C[x,y,z]$.

\vspace{3mm}

\begin{Defi}
Given a point $(k_1,k_2)\in\Z^2_{\geq 0}$ and a positive integer $\ell\in\Z_{>0}$, we define the polynomial
\[C(k_1,k_2,\ell)=\sum_{i_1\geq k_1,\ i_2\geq k_2,\ i_1+i_2\leq\ell}\binom{\ell}{i_1}\binom{\ell-i_1}{i_2}S(i_1,k_1)S(i_2,k_2)x^{i_1}y^{i_2}z^{\ell-i_1-i_2}\in\Z[x,y,z].\]
\label{polC}
\end{Defi}
Two direct properties about the polynomials $C(k_1,k_2,\ell)$ are
\begin{equation}
\begin{array}{l}
C(k_1,k_2,\ell)=0\mbox{ whenever }k_1+k_2>\ell,\\
\\
C(k_1,k_2,k_1+k_2)=\binom{k_1+k_2}{k_1}x^{k_1}y^{k_2}.\\
\end{array}
\label{eqC}
\end{equation}

The following alternative definition of the polynomials $C(k_1,k_2,\ell)$ can be deduced straightforwardly from the property (\ref{FormulaKT}) of Stirling numbers,
\begin{equation}
C(k_1,k_2,\ell)=\left\langle\frac{t^\ell}{\ell!}\right\rangle\frac{1}{k_1!k_2!}(e^{xt}-1)^{k_1}(e^{yt}-1)^{k_2}e^{zt},
\label{C2def}
\end{equation}
where $\left\langle\frac{t^\ell}{\ell!}\right\rangle f(t)$ denotes the coefficient of $\frac{t^\ell}{\ell!}$ in the expansion of the polynomial $f(t)$.

\vspace{2mm}

\begin{Defi}
Given a finite and ordered set of points $R\subseteq\Z^2_{\geq 0}$, consider the following matrix
\[M_R=\Big( C(k_1,k_2,\ell) \Big)_{(k_1,k_2)\in R,\ 0\leq\ell\leq r-1}\in\mathcal M_{r\times r}\big(\Z[x,y,z]\big),\]
where $r$ is the number of points in $R$.
\label{DefM}
\end{Defi}

Notice that each point $(k_1,k_2)\in R$ corresponds to the following row of the matrix
\[\Big(C(k_1,k_2,0),\ldots,C(k_1,k_2,r-1)\Big),\]
and the rows in the matrix $M_R$ are ordered as the points in $R$.

\vspace{3mm}

We are interested in the following linear system of polynomial equations
\begin{equation}
M_R\left(\begin{array}{c}
a_0\\
\vdots\\
a_{r-1}\\
\end{array}\right)=-\left(\begin{array}{c}
\vdots\\
C(k_1,k_2,r)\\
\vdots\\
\end{array}\right).
\label{SistStir}
\end{equation}
We will prove that, under certain conditions on the set $R$, the system has a solution ${\bf a}=(a_0,\ldots,a_{r-1})$ such that the polynomial $b_{\bf a}(s)=s^r+a_{r-1}s^{r-1}+\cdots+a_1s+a_0$ factors completely in $\C[x,y,z]$, and whose roots are somehow related to  the points in $R$. To state this relation precisely we need the following definition.

\begin{Defi}
Given $R$ a finite set of points in $\Z^2$, we define the monic polynomial
\[b_R(s)=\prod_{(i,j)\in R}(s-A_{ij})\]
where
\[A_{ij}=ix+jy+z\in\Z[x,y,z].\]
\label{DefbR}
\end{Defi}
We will prove in Theorem \ref{TheSt} that $b_{\bf a}(s)=b_R(s)$ (recall Notation \ref{Not2}). Or, in other words, we encode the solution to the system by giving the roots of the polynomial $b_{\bf a}(s)$. If we want the solution explicitly, all we have to do is use the elementary symmetric functions.

\vspace{3mm}

\begin{Exam}
The simplest example is the system obtained by imposing the conditions $x=1$ and $y=z=0$. Then, for $0\leq i\leq \ell$, the polynomials $C(i,0,\ell)$  reduce to
\[C(i,0,\ell)=S(\ell,i),\]
and then considering the system associated with the set $R=\{(0,0),(1,0),\ldots,(n-1,0)\}$, we recover the system in (\ref{SystemS1}). Moreover we have that $A_{ij}=i$ and we already know that $b_{\bf a}(s)=b_R(s)$ is a solution to the system.
\end{Exam}
In fact, as the example above suggests, we will consider compact sets of points $R\subset\Z^2_{\geq 0}$.

\begin{Defi}
Given a monomial ideal $J$, the set of standard monomials, denoted by ${\rm std}_J$, is the set of monomials which do not belong to $J$. We say that a finite set of points $R$ in $\Z^2$ satisfies the monomial condition if it can be identified with ${\rm std}_J$, for a zero-dimensional monomial ideal $J$ in $\C[x,y]$, by means of the identification
$$x^ay^b\longleftrightarrow (a,b).$$
\label{monomial}
\end{Defi}
Note that the fact that $J$ is zero-dimensional implies that ${\rm std}_J$ is finite.

\vspace{4mm}

Analogously as we defined the squared matrix $M_R$, we can extend Definition \ref{DefM} and define a matrix associated to the ordered set $R$, and of any length. For any $\ell\in\Z_{>0}$ we consider the matrix
\[M_{R,\ell}=\Big( C(k_1,k_2,j) \Big)_{(k_1,k_2)\in R,\ 0\leq j\leq\ell}\in\mathcal M_{r\times(\ell+1)}\Big(\Z[x,y,z]\Big).\]
With this notation we have that $M_R=M_{R,r-1}$, and we can write the system (\ref{SistStir}) in the equivalent form
\[M_{R,r}\left(\begin{array}{c}
a_0\\
\vdots\\
a_{r-1}\\
1\\
\end{array}\right)=\bar 0,\]
or in compact form
\[M_{R,r}({\bf a},1)^t=\overline 0,\]
where recall our notation ${\bf a}=(a_0,\ldots,a_{r-1})$. We write simply ${\bf a}$ instead of ${\bf a}_r$ to simplify notation, since $r$ is clear from the context.

\vspace{3mm}

Next we study the systems
\begin{equation}
M_{R,\ell}\big({\bf a},1\big)^t=\bar 0
\label{SystemGeneral}
\end{equation}
with $\ell\geq r$. Before solving these systems, we are going to illustrate how we encountered them. We summarize very briefly the process, see \cite{CC} for details.

When computing global $b$--functions we use the method of indeterminate coefficients to compute the minimal polynomial of the operator
\[s=\omega_1x_1\partial_1+\cdots+\omega_nx_n\partial_n\]
in $D/{\rm in}_{(-\omega,\omega)}(I)$, where $I$ is a holonomic ideal in the Weyl algebra $D$ and $\omega\in\R^n\setminus\{0\}$. We start computing powers of $s$ until there is a linear relation such that
\[s^r+a_{r-1}s^{r-1}+\cdots+a_0\in{\rm in}_{(-\omega,\omega)}(I)\]
for certain $a_{r-1},\ldots,a_0$. If the ideal $I$ is a hypergeometric ideal $H_A(\beta)$ with $A$ a matrix of size $n\times(n+2)$, we can reduce the problem to find a linear relation
\[\widetilde s^r+a_{r-1}\widetilde s^{r-1}+\cdots+a_0\in{\rm in}_{(-\omega,\omega)}(I)\]
where
\[\widetilde s=\alpha_ix_i\partial_i+\alpha_j\partial_j+\alpha_0\]
for certain $1\leq i<j\leq n$, and certain coefficients $\alpha_i,\alpha_j$ and $\alpha_0$, which are linear functions on $\beta$ and $\omega$.

Moreover there exists a set $B$ of monomials in $\C[\partial_i,\partial_j]$ contained in a Gröbner basis of ${\rm in}_{(-\omega,\omega)}(I)$ such that, for some $c_i,c_j>0$, $\partial_i^{c_i},\partial_j^{c_j}\in B$. Then it suffices to study when
\begin{equation}
\widetilde s^r+a_{r-1}\widetilde s^{r-1}+\cdots+a_0\equiv 0\mbox{ mod }B
\label{eqMallorca}
\end{equation}
Now, using that $\big(x\partial)^\ell=\sum_{k=0}^\ell S(\ell,k)x^k\partial^k$, we deduce that
\[\widetilde s^\ell=\sum_{k_1\geq 0,k_2\geq 0,k_1+k_2\leq\ell} C(k_1,k_2,\ell)x_i^{k_1}\partial_i^{k_1}x_j^{k_2}\partial_j^{k_2}\]
where the polynomials $C(k_1,k_2,\ell)$ are here specialized to the point $(x,y,z)=(\alpha_i,\alpha_j,\alpha_0)$. Solving equation (\ref{eqMallorca}) is equivalent to solving the system $M_{R,r}\big({\bf a},1\big)^t=0$.

\vspace{3mm}

\begin{The}
Let $R$ be a set of $r$ ordered points in $\Z^2$ satisfying the monomial condition. Consider the system
\[M_{R,\ell}\big({\bf a},1\big)^t=\bar 0\]
with $\ell\geq r$. Any ${\bf a}=(a_0,\ldots,a_{\ell-1})$ such that
\[b_{\bf a}(s)=\prod_{(i,j)\in R}(s-A_{ij})^{n_{ij}}\]
with $n_{ij}>0$ for any $(i,j)\in R$, and with $\sum_{(i,j)\in R} n_{ij}=\ell$, is a solution to the system.
\label{TheSt}
\end{The}

\begin{Rem}
Obviously $\big(b_{\bf a}(s)\big)_{red}=b_R(s)$, and in particular, when $\ell=r$, the solution ${\bf a}$  satisfies
\[b_{\bf a}(s)=b_R(s).\]
\end{Rem}

\vspace{3mm}

Before going further we need some combinatorial results regarding the polynomials $C(k_1,k_2,\ell)$. First let us recall the {\em Vandermonde Convolution} formula:
\begin{equation}
\binom{n+m}{k}=\sum_{j=0}^k\binom{n}{j}\binom{m}{k-j}.
\label{VDMD}
\end{equation}

\vspace{2mm}

\begin{Defi}
Given $(k_1,k_2)\in\N^2$ we define the sets
\[\begin{array}{l}
S_{k_1,k_2}=\{(i,j)\ |\ 0\leq i\leq k_1, 0\leq j\leq k_2\}\\
\\
S_{k_1,k_2}^\bullet=S_{k_1,k_2}\setminus\{(k_1,k_2)\}\\
\end{array}\]
Moreover we denote
\[\begin{array}{l}
c_{ij}=\binom{k_1}{i}\binom{k_2}{j}(-1)^{k_1+k_2-i-j}\\
\\
d_{ij}=\binom{k_1}{i}\binom{k_2}{j}i!j!\\
\end{array}\]
Notice that the notations $c_{ij}$ and $d_{ij}$ depend on $(k_1,k_2)$, but in order to simplify notation, we do not write this dependence explicitly.
\end{Defi}

\vspace{3mm}

Next we prove a kind of inversion formulas among the polynomials $C(k_1,k_2,\ell)$ and $A_{ij}$.

\begin{Lem}
For all $(k_1,k_2)\in\Z^2_{\geq 0}$ and $\ell\geq 0$ we have
\begin{enumerate}
\item[(i)] $\sum_{(i,j)\in S_{k_1,k_2}}d_{ij}C(i,j,\ell)=A_{k_1k_2}^\ell$.

\

\item[(ii)] $k_1!k_2!C(k_1,k_2,\ell)=\sum_{(i,j)\in S_{k_1,k_2}}c_{ij}A_{ij}^\ell$.

\

\item[(iii)] $\sum_{(i,j)\in S_{k_1,k_2}^\bullet}c_{ij}A_{ij}^\ell=-\sum_{(i,j)\in S_{k_1,k_2}^\bullet}d_{ij}C(i,j,\ell)$.
\end{enumerate}
\label{combinatorics}
\end{Lem}

{\em Proof.} We prove the identity in (i). We have that
\[\sum_{(i,j)\in S_{k_1,k_2}}\binom{k_1}{i}\binom{k_2}{j}i!j!C(i,j,\ell)=\]
\[\begin{array}{l}
=\sum_{(i,j)\in S_{k_1,k_2}}\binom{k_1}{i}\binom{k_2}{j}i!j!\sum_{i_1\geq i,i_2\geq j,i_1+i_2\leq \ell}\binom{\ell}{i_1}\binom{\ell-i_1}{i_2}S(i_1,i)S(i_2,j)x^{i_1}y^{i_2}z^{\ell-i_1-i_2}\\
\\
=\sum_{(i_1,i_2)\in T_\ell}\binom{\ell}{i_1}\binom{\ell-i_1}{i_2}x^{i_1}y^{i_2}z^{\ell-i_1-i_2}\sum_{(i,j)\in S_{i_1,i_2}\cap S_{k_1,k_2}}\binom{k_1}{i}\binom{k_2}{j}i!j!S(i_1,i)S(i_2,j)\\
\end{array}\]
where $T_\ell=\{(i,j)\in\Z^2_{\geq 0}\ |\ i+j\leq \ell\}$. Now notice that
\[\begin{array}{ll}
\sum_{(i,j)\in S_{i_1,i_2}\cap S_{k_1,k_2}}\binom{k_1}{i}\binom{k_2}{j}i!j!S(i_1,i)S(i_2,j) & =\sum_{i=0}^{min\{i_1,k_1\}}\sum_{j=0}^{min\{i_2,k_2\}}\binom{k_1}{i}\binom{k_2}{j}i!j!S(i_1,i)S(i_2,j)\\
 \\
 & =\sum_{i=0}^{i_1}\sum_{j=0}^{i_2}\binom{k_1}{i}\binom{k_2}{j}i!j!S(i_1,i)S(i_2,j)\\
 \\
 & =\left(\sum_{i=0}^{i_1}\binom{k_1}{i}i!S(i_1,i)\right)\left(\sum_{j=0}^{i_2}\binom{k_2}{j}j!S(i_2,j)\right)\\
 \\
 & = k_1^{i_1}k_2^{i_2}\\
 \end{array}\]
where in the last equality we use the identity:
\[x^n=\sum_{i=0}^n\ S(n,i)x(x-1)\cdots (x-i+1)\]
Therefore we are done since
\[\sum_{(i,j)\in S_{k_1,k_2}}\binom{k_1}{i}\binom{k_2}{j}i!j!C(i,j,\ell)=\sum_{(i_1,i_2)\in T_\ell}\binom{\ell}{i_1}\binom{\ell-i_1}{i_2}x^{i_1}y^{i_2}z^{\ell-i_1-i_2}k_1^{i_1}k_2^{i_2}=(k_1x+k_2y+z)^\ell.\]

\vspace{3mm}

Now we prove the identity in (ii). We have that
\[\sum_{i=0}^{k_1}\sum_{j=0}^{k_2}\binom{k_1}{i}\binom{k_2}{j}(-1)^{k_1+k_2-i-j}(ix+jy+z)^\ell=\]
\[=\sum_{i=0}^{k_1}\sum_{j=0}^{k_2}\sum_{(a,b)\in T_\ell}\binom{k_1}{i}\binom{k_2}{j}(-1)^{k_1+k_2-i-j}\binom{\ell}{a}\binom{\ell-a}{b}x^ay^bz^{\ell-a-b}i^aj^b=\]
\[=\sum_{(a,b)\in T_\ell}\binom{\ell}{a}\binom{\ell-a}{b}x^ay^bz^{\ell-a-b}\left(\sum_{i=0}^{k_1}\binom{k_1}{i}(-1)^{k_1-i}i^a\right)\left(\sum_{j=0}^{k_2}\binom{k_2}{j}
(-1)^{k_2-j}j^b\right).\]
We use the following closed form of the Stirling numbers of second kind:
\[S(n,m)=\frac{1}{m!}\sum_{k=0}^m(-1)^{m-k}\binom{m}{k}k^n\]
and, since $S(n,m)=0$ when $n<m$, we are done.

\vspace{3mm}

The identity in (iii) follows directly by (i) and (ii). \hfill$\Box$

\vspace{6mm}

{\em Proof of Theorem \ref{TheSt}.}
 The proof is by double induction on $r>0$ and $\ell\geq r$.

For $r=1$, the set of points is $R=\{(0,0)\}$ and hence the system under study is simply the equation
\[\Big(C(0,0,0)\ C(0,0,1)\ \cdots\ C(0,0,\ell)\Big)\left(\begin{array}{c}
a_0\\
\vdots\\
a_{\ell-1}\\
1\\
\end{array}\right)=0,\]
that is, the equation
\[a_0+a_1z+\cdots+a_{\ell-1}z^{\ell-1}+z^\ell=0.\]
It is obvious that
\[a_j=\binom{\ell}{j}(-1)^{\ell-j}z^{\ell-j},\ \ 0\leq j<\ell,\]
which are the coefficients of the polynomial $(s-z)^\ell=(s-A_{00})^\ell$, gives a solution of the previous equation.

\vspace{2mm}

For $r>1$, suppose first that $\ell=r$. We have to prove that ${\bf a}=(a_0,\ldots,a_{r-1})$ with $b_{\bf a}(s)=b_R(s)$ is a solution to the system $M_{R,r}\big({\bf a},1\big)^t=\bar 0$. We can always choose $(k_1,k_2)\in R$ such that the set $R':=R\setminus\{(k_1,k_2)\}$ has the monomial condition. Then, by induction hypothesis we have that $\mbox{\boldmath{$\l$}}=(\l_0,\ldots,\l_{r-2})$ such that $b_{\mbox{\boldmath{$\l$}}}(s)=b_{R'}(s)$, is a solution to the system
\[M_{R',r-1}\big(\mbox{\boldmath{$\l$}},1\big)^t=\bar 0.\]
Notice that $b_R(s)=(s-A_{k_1k_2})b_{R'}(s)$, and hence we have to prove that
\begin{equation}
M_{R,r}\left[\left(\begin{array}{c}
0\\
\l_0\\
\vdots\\
\l_{r-3}\\
\l_{r-2}\\
1\\\end{array}\right)-A_{k_1k_2}\left(\begin{array}{c}
\l_0\\
\l_1\\
\vdots\\
\l_{r-2}\\
1\\
0\\
\end{array}\right)\right]=\bar 0.
\label{eqC'}
\end{equation}
We can write the matrix $M_{R,r}$ as
\[M_{R,r}=\left(\begin{array}{ccc}
M_{R',r-1} &  & \begin{array}{c}
           \vdots\\
           C(\star,r)\\
           \vdots\\
           \end{array}\\
\\
0\ \cdots\ 0\ C(k_1,k_2,k_1+k_2) & \cdots & C(k_1,k_2,r)
\end{array}\right)\]
where $\star$ runs over the points in $R'$. Let us consider the system $M_{R',r}\big({\bf a},1\big)^t=\bar 0$. We know by induction hypothesis that for any $(i,j)\in R'$, the coefficients of $(s-A_{ij})b_{R'}(s)$ give a solution to this system. Then
\[\left(\begin{array}{cccc}
 & & & \vdots\\
 \\
 & M_{R',r-1} & & C(\star,r)\\
 \\
 & & & \vdots\\
 \end{array}\right)
\left[\left(\begin{array}{c}
0\\
\lambda_0\\
\vdots\\
\lambda_{r-2}\\
1\\
\end{array}\right)-A_{ij}\left(\begin{array}{c}
\lambda_0\\
\vdots\\
\lambda_{r-2}\\
1\\
0\\
\end{array}\right)\right]=\bar 0.\]
Since $M_{R',r}\big(\l_0,\ldots,\l_{r-2},1,0\big)^t=\bar 0$ is equivalent to $M_{R',r-1}\big(\l_0,\ldots,\l_{r-2},1\big)^t=\bar 0$, we deduce that
\begin{equation}
\left(\begin{array}{cccc}
 & & & \vdots\\
 \\
 & M_{R',r-1} & & C(\star,r)\\
 \\
 & & & \vdots\\
 \end{array}\right)
 \left(\begin{array}{c}
 0\\
 \l_0\\
 \vdots\\
 \l_{r-2}\\
 1\\
 \end{array}\right)=\bar 0.
 \label{eqs}
 \end{equation}
Then, to solve system (\ref{eqC'}), we only have to prove the equation corresponding to the point $(k_1,k_2)$:
\begin{equation}
\big(0,\ldots,0,C(k_1,k_2,k_1+k_2),\ldots,C(k_1,k_2,r)\big)\left[\left(\begin{array}{c}
0\\
\l_0\\
\vdots\\
\l_{r-2}\\
1\\
\end{array}\right)-A_{k_1k_2}\left(\begin{array}{c}
\l_0\\
\l_1\\
\vdots\\
1\\
0\\
\end{array}\right)\right]=\bar 0,
\label{eqfinal}
\end{equation}
or in other words, we have to prove the following equality
\[\sum_{\ell=k_1+k_2}^{r}C(k_1,k_2,\ell)\l_{\ell-1}=A_{k_1k_2}\sum_{\ell=k_1+k_2}^{r}C(k_1,k_2,\ell)\l_\ell,\]
where $\l_r:=0$ and recall that $\l_{r-1}=1$.
By Lemma \ref{combinatorics} (ii), this equality is
\[\sum_{\ell=k_1+k_2}^{r}\sum_{(i,j)\in S_{k_1,k_2}}c_{ij}A_{ij}^\ell \l_{\ell-1}=A_{k_1k_2}\sum_{\ell=k_1+k_2}^{r}\sum_{(i,j)\in S_{k_1,k_2}}c_{ij}A_{ij}^\ell \l_\ell.\]
 We write our equation as follows
\[\sum_{\ell=k_1+k_2}^{r}c_{k_1k_2}A_{k_1k_2}^\ell \l_{\ell-1}+\sum_{\ell=k_1+k_2}^{r}\sum_{(i,j)\in S_{k_1,k_2}^\bullet}c_{ij}A_{ij}^\ell \l_{\ell-1}=\sum_{\ell=k_1+k_2}^{r-1}c_{k_1k_2}A_{k_1k_2}^{\ell+1}\l_\ell+A_{k_1k_2}\sum_{\ell=k_1+k_2}^{r-1}\sum_{(i,j)\in S_{k_1,k_2}^\bullet}c_{ij}A_{ij}^\ell \l_\ell.\]
After some cancelations the equation looks like
\begin{equation}\label{sinnombre}
c_{k_1k_2}A_{k_1k_2}^{k_1+k_2}\l_{k_1+k_2-1}+\sum_{\ell=k_1+k_2}^{r}\sum_{(i,j)\in S_{k_1,k_2}^\bullet}c_{ij}A_{ij}^\ell \l_{\ell-1}=A_{k_1k_2}\sum_{\ell=k_1+k_2}^{r-1}\sum_{(i,j)\in S_{k_1,k_2}^\bullet}c_{ij}A_{ij}^\ell \l_\ell.
\end{equation}
By Lemma \ref{combinatorics} (i) and (iii) we have
\[\sum_{(i,j)\in S_{k_1,k_2}^\bullet}c_{ij}A_{ij}^\ell=-\left(A_{k_1k_2}^\ell-k_1!k_2!C(k_1,k_2,\ell)\right),\]
and applying it for $\ell=k_1+k_2-1$ we deduce that
\[\sum_{(i,j)\in S_{k_1,k_2}}c_{ij}A_{ij}^{k_1+k_2-1}=0.\]
Then $c_{k_1k_2}A_{k_1k_2}^{k_1+k_2}=-A_{k_1k_2}\sum_{(i,j)\in S_{k_1,k_2}^\bullet}c_{ij}A_{ij}^{k_1+k_2-1}$, and we can write equation (\ref{sinnombre}) as
\begin{equation}
\sum_{\ell=k_1+k_2}^{r}\left(\sum_{(i,j)\in S_{k_1,k_2}^\bullet}c_{ij}A_{ij}^\ell\right)\l_{\ell-1}=A_{k_1k_2}\sum_{\ell=k_1+k_2-1}^{r-1}\left(\sum_{(i,j)\in S_{k_1,k_2}^\bullet}c_{ij}A_{ij}^\ell\right)\l_\ell.
\label{eqQuasifin}
\end{equation}
Notice that we have reduced our problem to an equation in terms of $S_{k_1,k_2}^\bullet$, which is contained in $R'$, the set where we can apply induction. Let us write now equation (\ref{eqQuasifin}) again in terms of $C(i,j,\ell)$, by using Lemma \ref{combinatorics} (iii). We have
\begin{equation}\label{july2}
\sum_{(i,j)\in S_{k_1,k_2}^\bullet}d_{ij}\left(\sum_{\ell=k_1+k_2}^rC(i,j,\ell)\l_{\ell-1}\right)
=A_{k_1k_2}\sum_{(i,j)\in S_{k_1,k_2}^\bullet}d_{ij}\left(\sum_{\ell=k_1+k_2-1}^{r-1}C(i,j,\ell)\l_\ell\right).
\end{equation}
Since $S_{k_1,k_2}^\bullet\subseteq R'$, from (\ref{eqs}) we have that for any $(i,j)\in S_{k_1,k_2}^\bullet\subseteq R'$
\[\big(0,\ldots,0,C(i,j,i+j),\ldots,C(i,j,r)\big)\left(\begin{array}{c}
0\\ \l_0\\ \vdots\\ \l_{r-2}\\ 1\\ \end{array}\right)=0,\]
and by induction hypothesis, for any $(i,j)\in S_{k_1,k_2}^\bullet\subseteq R'$ we have that
\[\big(0,\ldots,0,C(i,j,i+j),\ldots,C(i,j,r-1)\big)\left(\begin{array}{c}
\l_0\\ \vdots\\ \l_{r-2}\\ 1\\ \end{array}\right)=0.\]
We deduce that (\ref{july2}) is equivalent to
\begin{equation}\label{july}
\sum_{(i,j)\in S_{k_1,k_2}^\bullet}\sum_{\ell=1}^{k_1+k_2-1}d_{ij}C(i,j,\ell)\l_{\ell-1}=A_{k_1k_2}\sum_{(i,j)\in S_{k_1,k_2}^\bullet}\sum_{\ell=0}^{k_1+k_2-2}d_{ij}C(i,j,\ell)\l_\ell
\end{equation}
%where we set $\l_{-1}=0$,
since for any $(i,j)\in S_{k_1,k_2}^\bullet$ we have $i+j<k_1+k_2$.
Then equation (\ref{july}) can be written as
\[\sum_{\ell=1}^{k_1+k_2-1}\left(\sum_{(i,j)\in S_{k_1,k_2}^\bullet}d_{ij}C(i,j,\ell)\right)\l_{\ell-1}=A_{k_1k_2}\sum_{\ell=0}^{k_1+k_2-2}\left(\sum_{(i,j)\in S_{k_1,k_2}^\bullet}d_{ij}C(i,j,\ell)\right)\l_\ell.\]
Now we claim that for $1\leq \ell\leq k_1+k_2-1$
\[\sum_{(i,j)\in S_{k_1,k_2}^\bullet}d_{ij}C(i,j,\ell)=A_{k_1k_2}\sum_{(i,j)\in S_{k_1,k_2}^\bullet}d_{ij}C(i,j,\ell-1)\]
and hence the equation we wanted to prove holds.

\vspace{2mm}

The claim follows by Lemma \ref{combinatorics} (i), since we have
\[\begin{array}{ll}
\sum_{(i,j)\in S_{k_1,k_2}^\bullet}d_{ij}C(i,j,\ell) & =A_{k_1k_2}^\ell-k_1!k_2!C(k_1,k_2,\ell) =A_{k_1k_2}^\ell\\
\\
 &=A_{k_1k_2}(A_{k_1k_2}^{\ell-1}-k_1!k_2!C(k_1,k_2,\ell-1))\\
\\
 & =A_{k_1k_2}\sum_{(i,j)\in S_{k_1,k_2}^\bullet}d_{ij}C(i,j,\ell-1)\\
 \end{array}\]
where we used that $C(k_1,k_2,\ell)=0$ if $\ell<k_1+k_2$.

\vspace{5mm}

Suppose now that $\ell>r$. We have to prove that the system
\[M_{R,\ell}\big({\bf a},1\big)^t=\bar 0\]
has a solution ${\bf a}$ such that $b_{\bf a}(s)=\prod_{(i,j)\in R}(s-A_{ij})^{n_{ij}}$ with $\sum n_{ij}=\ell$. By induction hypothesis we have that any $\mbox{\boldmath{$\l$}}=(\l_0,\ldots,\l_{\ell-2},1)$ such that $b_{\mbox{\boldmath{$\l$}}}(s)=\prod_{(i,j)\in R}(s-A_{ij})^{m_{ij}}$ with $\sum m_{ij}=\ell-1$ satisfies
\begin{equation}
M_{R,\ell-1}\left(\begin{array}{c}
\l_0\\
\vdots\\
\l_{\ell-2}\\
1\\
\end{array}\right)=\bar 0.
\label{eqJueves}
\end{equation}
 We have to prove that for any $(i,j)\in R$ the coefficients of the polynomial
\[(s-A_{ij})(s^{\ell-1}+\l_{\ell-2}s^{\ell-2}+\cdots+\l_0)\]
give a solution to the system $M_{R,\ell}\big({\bf a},1\big)^t=\bar 0$. Then we have to prove that
\[M_{R,\ell}\left[\left(\begin{array}{c}
0\\
\l_0\\
\vdots\\
\l_{\ell-2}\\
1\\
\end{array}\right)-A_{ij}\left(\begin{array}{c}
\l_0\\
\vdots\\
\l_{\ell-2}\\
1\\
0\\
\end{array}
\right)\right]=\bar 0.\]
By (\ref{eqJueves}) it is clear that $-A_{ij}M_{R,\ell}\big(\l_0,\ldots,\l_{\ell-2},1,0\big)^t=\bar 0$ and hence we have to prove
\[M_{R,\ell}\big(0,\l_0,\ldots,\l_{\ell-2},1\big)^t=\bar 0.\]
We are going to prove that, for any $(k_1,k_2)\in R$, we have
\[\sum_{n=1}^\ell C(k_1,k_2,n)\l_{n-1}=0,\]
with the convention $\l_{\ell-1}=1$.
By Lemma \ref{combinatorics} (ii) this equation can be written as
\[\sum_{n=1}^\ell\frac{1}{k_1!k_2!}\sum_{(i,j)\in S_{k_1,k_2}^\bullet}c_{ij}A_{ij}^n\l_{n-1}=0.\]
Or equivalently
\[\sum_{(i,j)\in S_{k_1,k_2}^\bullet}c_{ij}A_{ij}\sum_{n=1}^\ell A_{ij}^{n-1}\l_{n-1}=0.\]
This last equality holds, since $\sum_{n=1}^\ell A_{ij}^{n-1}\l_{n-1}=b_{\mbox{\boldmath{$\l$}}}(A_{ij})$ and this vanishes since $(i,j)\in R$.
\hfill$\Box$

\vspace{3mm}

\begin{Rem}
It can be checked that  the condition $n_{ij}>0$ is necessary in the statement of Theorem \ref{TheSt}.
\end{Rem}

\vspace{3mm}

Next we study the determinant of the squared matrix $M_R$ (recall that $M_R=M_{R,r-1}$).
The determinant of $M_R$ is a polynomial in $\Z[x,y,z]$. We prove that the polynomial is in fact in $\Z[x,y]$ and that it is not identically zero.

\begin{Pro}
Let $R$ be a set of ordered points in $\Z^2$ satisfying the monomial condition and with cardinal at least two. Then
\[det(M_R)=\prod_{(i,j)\in R}\frac{1}{i!j!}\cdot\prod_{(i,j)<(i',j')\in R}(A_{ij}-A_{i'j'})\]
where $<$ orders the points in $R$ as the corresponding rows of the matrix $M_R$ are ordered.
\label{determinant}
\end{Pro}

{\em Proof .} (This proof was provided to the authors by Christian Krattenthaler).

\vspace{2mm}

By the alternative definition of $C(k_1,k_2,\ell)$ given in (\ref{C2def}) we have that
\begin{equation}
\mbox{det}(M_R)=\prod_{(k_1,k_2)\in R}\frac{1}{k_1!k_2!}\cdot\mbox{ det}\left(\left\langle\frac{t^\ell}{\ell!}\right\rangle(e^{xt}-1)^{k_1}(e^{yt}-1)^{k_2}e^{zt}\right)_{(k_1,k_2)\in R,\ 0\leq \ell<r}.
\label{eqK1}
\end{equation}
We claim that the determinant is equal to
\begin{equation}
\prod_{(k_1,k_2)\in R}\frac{1}{k_1!k_2!}\cdot\mbox{ det}\left(\left\langle\frac{t^\ell}{\ell!}\right\rangle e^{(k_1x+k_2y+z)t}\right)_{(k_1,k_2)\in R,\ 0\leq \ell<r}.
\label{eqK2}
\end{equation}
Then we are done, because this determinant can be written as
\[\prod_{(k_1,k_2)\in R}\frac{1}{k_1!k_2!}\cdot\mbox{ det}\Big((k_1x+k_2y+z)^\ell\Big)_{(k_1,k_2)\in R,\ 0\leq \ell<r}.\]
This is a Vandermonde determinant, which can therefore be evaluated, and the result follows.

\vspace{2mm}

Now we prove the claim. We impose the lexicographic order on the points in $R$, and denote it by $\preceq$.
Fix a point $(m_1,m_2)\in R$. We are going to prove that the determinant in the right hand side of the identity in (\ref{eqK1}) equals
\begin{equation}
\mbox{det}\left(\begin{array}{ll}
\left\langle\frac{t^\ell}{\ell!}\right\rangle e^{(k_1x+k_2y+z)t} & \mbox{ for }(k_1,k_2)\preceq (m_1,m_2)\\
\\
\left\langle\frac{t^\ell}{\ell!}\right\rangle (e^{xt}-1)^{k_1}(e^{yt}-1)^{k_2}e^{zt} & \mbox{ for }(k_1,k_2)\not\preceq(m_1,m_2)\\
\end{array}\right)_{(k_1,k_2)\in R,\ 0\leq \ell<r}.
\label{eqK4}
\end{equation}
We prove this by induction on the elements in $R$ with respect to the order $\preceq$. Clearly, the claim is correct for $(m_1,m_2)=(0,0)$, which gives the start of the induction. For the induction step, let $(M_1,M_2)$ denote the successor of $(m_1,m_2)$ in $R$, with respect to the (total) order $\preceq$. Our determinant looks, by induction hypothesis,
\[\mbox{det}\left(\begin{array}{ll}
\left\langle\frac{t^\ell}{\ell!}\right\rangle e^{(k_1x+k_2y+z)t} & \mbox{ for }(k_1,k_2)\preceq (m_1,m_2)\\
\\
\left\langle\frac{t^\ell}{\ell!}\right\rangle (e^{xt}-1)^{M_1}(e^{yt}-1)^{M_2}e^{zt} &\\
\\
\left\langle\frac{t^\ell}{\ell!}\right\rangle (e^{xt}-1)^{k_1}(e^{yt}-1)^{k_2}e^{zt} & \mbox{ for }(k_1,k_2)\not\preceq(M_1,M_2)\\
\end{array}\right)_{(k_1,k_2)\in R,\ 0\leq \ell<r}.\]
 We have
\begin{equation}
\begin{array}{ll}
(e^{xt}-1)^{M_1}(e^{yt}-1)^{M_2}e^{zt} & =\sum_{k_1\leq M_1,k_2\leq M_2}(-1)^{M_1+M_2-k_1-k_2}\binom{M_1}{k_1}\binom{M_2}{k_2}e^{(k_1x+k_2y+z)t}\\
\\
 & =e^{(M_1x+M_2y+z)t}+\sum_{\stackrel{k_1\leq M_1,k_2\leq M_2}{(k_1,k_2)\neq(M_1,M_2)}}(-1)^{M_1+M_2-k_1-k_2}\binom{M_1}{k_1}\binom{M_2}{k_2}e^{(k_1x+k_2y+z)t}\\
 \end{array}
\label{eqK3}
\end{equation}
We may now substract appropriate linear combinations of the rows of the matrix in (\ref{eqK4}) indexed by $(k_1,k_2)\preceq(m_1,m_2)$ from the row indexed by $(M_1,M_2)$ to eliminate all terms $e^{(k_1x+k_2y+z)t}$ in the expansion (\ref{eqK3}) with $(k_1,k_2)\preceq(m_1,m_2)$. In other words, all terms in (\ref{eqK3}) are eliminated except for the leading term $e^{(M_1x+M_2y+z)t}$. Thus we have established our claim for $(M_1,M_2)$, the successor of $(m_1,m_2)$.

 This finishes the proof of the equality among (\ref{eqK1}) and (\ref{eqK2}). We just have to choose the (lexicographically) largest element of $R$ for $(m_1,m_2)$ in (\ref{eqK4}).
\hfill$\Box$

\vspace{4mm}

We end the section illustrating why the condition we impose on the set $R$ is necessary.

\begin{Rem}
Notice that if we do not ask the monomial condition of the set of points $R$, the result in Proposition \ref{determinant} is not true any longer. Indeed, consider the set
\[R=\{(0,0),(1,0),(2,0),(3,0),(0,1),(0,2),(0,4)\}.\]
It does not have the monomial condition since $(0,3)\notin R$. We can check that
\[\mbox{det}(M_R)=-2x^6y^7(2x-y)(3x-y)(x-2y)(3x-2y)(x-y)^2(11x^2-42xy+37y^2),\]
and $11x^2-42xy+37y^2$ is irreducible in $\Z[x,y]$.

\vspace{2mm}

This example also shows that Theorem \ref{TheSt} is not true if we drop the monomial condition. Indeed, the system $M_{R,7}\big({\bf a},1\big)^t=\bar 0$ is not solvable in $\C[x,y,z]$. It is, however, in $\C\big(x,y,z\big)$ and the function $b_{\bf a}(s)$ is
\[b_{\bf a}(s)=(s-z)(s-z-y)(s-z-2y)(s-z-3x)(s-z-x)(s-z-2x)\frac{p(s)}{11x^2-42xy+37y^2}\]
where $p(s)=11sx^2-42sxy+37sy^2+6x^3-175y^3-11x^2z-77x^2y+42xyz-37y^2z+222xy^2$.

\vspace{2mm}

Finally observe how the biggest subset of $R$ satisfying the monomial condition,
$$R_0=\{(0,0),(1,0),(2,0),(3,0),(0,1),(0,2)\},$$
appears both in the determinant and in the function $b_{\bf a}(s)$. More precisely,
\[det(M_R)=det(M_{R_0})\cdot y^4(11x^2-42xy+37y^2)\]
and
\[\begin{array}{ll}
b_{\bf a}(s) & =b_{R_0}(s)\cdot\frac{11sx^2-42sxy+37sy^2+6x^3-175y^3-11x^2z-77x^2y+42xyz-37y^2z+222xy^2}{11x^2-42xy+37y^2}\\
\\
 & =b_{R_0}(s)\Big(s-z+\frac{(3x-7y)(2x^2-21xy+25y^2)}{11x^2-42xy+37y^2}\Big)\\
 \end{array}\]
\end{Rem}

\vspace{3mm}

\section{Study of the system in some quotient rings of $\C[x,y,z]$}

If we consider systems of the form (\ref{SystemGeneral}) under a specialization of $x,y$ and $z$, we obtain a linear system of equations with complex coefficients. Along this section we will anyway denote by $M_{R,\ell}$, $A_{ij}$ and $b_R(s)$ the corresponding objects under the specialization considered. Obviously we get solutions to this system by simply specializing  the solutions obtained in Theorem \ref{TheSt}.

Recall that we explained briefly how these systems appeared when computing the global $b$-function of a hypergeometric ideal. More concretely when solving equations of the form
\[s^r+a_{r-1}s^{r-1}+\cdots+a_0=0,\]
with $r$ the smallest possible positive integer. In Section \ref{S2} we have solved it for $r$ equal to the cardinal of the set $R$ (which is related to the ideal ${\rm in}_{(-\omega,\omega)}(H_A(\beta))$ in a way too technical to make more explicit here). The question is, can we solve an analogous equation of smaller degree? This is the underlying motivation of this section.

\vspace{2mm}

We know solutions to the systems for $\ell\geq r$, where $r$ is the number of points of $R$, and the question is wether we can find solutions for $\ell<r$. We will find solutions to the systems of $r$ equations in $\ell<r$ variables for certain specializations, more concretely those such that det$(M_R)=0$.

\vspace{2mm}

By Proposition \ref{determinant} we have that  det$(M_R)=0$ is equivalent to $A_{ij}=A_{i'j'}$ for certain different points $(i,j),(i',j')\in R$. Let us then treat the problem with a bit more generality: instead of systems of linear equations with complex coefficients, we consider the systems of the form (\ref{SystemGeneral}) with polynomial coefficients as in previous section, but in the quotient ring $\C[x,y,x]/(A_{ij}-A_{i'j'})$ for certain $(i,j),(i',j')\in R$.

\vspace{2mm}

Note that det$(M_R)=0$ is equivalent to $b_{R}(s)$ non-reduced. Then, for $\ell\geq r=|R|$ the solutions ${\bf a}$ we find in Theorem \ref{TheSt} (under the specialization of the type $A_{ij}=A_{i'j'}$) correspond all to non-reduced polynomials $b_{\bf a}(s)$. We will see how, while looking for solutions to the system $M_{R,\ell}\big({\bf a},1\big)^t=\bar 0$ in $\C[x,y,z]/(A_{ij}-A_{i'j'})$ for $r<\ell$, we will answer the following related question: Can we find a reduced solution to the system for $\ell<r$? The answer is yes, as we will see below.

\vspace{2mm}

The condition $A_{ij}-A_{i'j'}=0$ can be written as:
\begin{equation}
ax+by=0,
\label{eqD0}
\end{equation}
with $a\in\Z_{\geq 0}$ and $b\in\Z$ and $(a,b)\neq (0,0)$.
This condition can be decomposed in the following cases:
\vspace{2mm}
\begin{enumerate}
\item[(i)] $ax+by=0$ with $a,b>0$.
\\
\item[(ii)] $ax+by=0$ with $a>0$ and $b<0$.
\\
\item[(iii)] $x=0$.
\\
\item[(iv)] $y=0$.
\end{enumerate}
\vspace{2mm}

From now on we will consider the systems $M_{R,\ell}\big({\bf a},1\big)^t=\bar 0$ in the ring $\C[x,y,z]/(ax+by)$ with $(a,b)\in\R^2\setminus\{(0,0)\}$.

\vspace{2mm}

The aim of this Section is to prove the following:

\begin{Pro}
Let $R$ be a set in $\Z^2$ satisfying the monomial condition, and consider the system
\[M_{R,\ell}\big({\bf a},1\big)^t=\bar 0\ \mbox{ in the ring }\ \C[x,y,z]/(ax+by),\]
for $\ell\geq r_0$, where $r_0$ is the degree of $\big(b_R(s)\big)_{red}$. There exists a set $R_{a,b}\subseteq R$ with cardinal $r_0$ and such that any ${\bf a}=(a_0,\ldots,a_{\ell-1})$ with
\[b_{\bf a}(s)=\prod_{(i,j)\in R_{a,b}}\big(s-A_{ij}\big)^{n_{ij}}\]
with $n_{ij}>0$ and $\sum n_{ij}=\ell$, is a solution to the system. Moreover, if $\ell=r_0$ the solution described above is
\[b_{\bf a}(s)=b_{R_{a,b}}(s)=\big(b_R(s)\big)_{red}\]
\label{PropSec3}
\end{Pro}

In general the set $R_{a,b}$ will not inherit the monomial condition from $R$, and therefore we can not use directly the results of Section \ref{S2}.
However we have the same kind of result as Proposition \ref{determinant} regarding the determinant of the matrix $M_{R_{a,b}}$.

\begin{Pro}
Let $R$ be a set of points in $\Z^2$ with the monomial condition. For any $\ell\in\Z_{>0}$ the systems $M_{R,\ell}\big({\bf a},1\big)^t=\bar 0$ and $M_{R_{a,b},\ell}\big({\bf a},1\big)^t=\bar 0$ are equivalent in the ring $\C[x,y,z]/(ax+by)$. Moreover we have
\[{\rm det}\big(M_{R_{a,b}}\big)=\prod_{(k_1,k_2)\in R_{a,b}}\frac{1}{k_1!k_2!}\prod_{(i,j)<(i'j')\in R_{a,b}}\big(A_{ij}-A_{i'j'}\big)\]
where $<$ orders the points in $R_{a,b}$ as they are ordered in the matrix (recall that, by definition, any point in $R_{a,b}$ corresponds to a row in $M_{R_{a,b}}$). In particular
\[{\rm det}\big(M_{R_{a,b}}\big)\neq 0\ \mbox{ in }\C[x,y,z]/(ax+by).\]
\label{Result2}
\end{Pro}

\vspace{2mm}

Notice that Lemma \ref{combinatorics} also works in $\C[x,y,z]/(ax+by)$, by simply imposing the condition $ax+by=0$ on the identities.

\vspace{3mm}

Next Proposition treats the cases (iii) and (iv), which are very simple to deal with.

\begin{Pro}
Let $R\subseteq\Z^2$ be a set of points with the monomial condition. In the ring $\C[x,y,z]/(x)$ (resp. in $\C[x,y,z]/(y)$), for any $\ell\in\Z_{\geq 0}$ the systems $M_{R,\ell}\big({\bf a},1\big)^t=\bar 0$ and $M_{R_0,\ell}\big({\bf a},1\big)^t=\bar 0$ are equivalent, where $R_0=\{(0,j)\in R\}$ (resp. $R_0=\{(i,0)\in R\}$). For any $\ell\geq r_0={\rm deg}\big(b_R(s)_{red}\big)$, any ${\bf a}$ such that $b_{\bf a}(s)=\prod_{(i,j)\in R_0}\big(s-A_{ij}\big)^{n_{ij}}$, with $n_{ij}>0$ and $\sum n_{ij}=\ell$ is a solution to $M_{R_0,\ell}\big({\bf a},1\big)^t=\bar 0$, and for $\ell=r_0$
\[b_{\bf a}(s)=b_{R_0}(s)=\big(b_R(s)\big)_{red}.\]
Moreover det$\big(M_{R_0}\big)=\prod_{j=0}^{r_0-1}y^j$.
\label{FACIL}
\end{Pro}

{\em Proof.}
Considering the system $M_{R,\ell}\big({\bf a},1\big)^t=\bar 0$ in $\C[x,y,z]/(x)$ (the case $\C[x,y,z]/(y)$ is completely analogous), we have that $A_{ij}=jy+z$ and therefore it is clear that
\[\big(b_R(s)\big)_{red}=b_{R_0}(s).\]
Hence the cardinal of $R_0$ equals $r_0$, the degree of $\big(b_R(s)\big)_{red}$. The equivalence of the two systems follows, since, for any $\ell\geq 0$:
\vspace{2mm}

$\bullet$ if $x=0$ then $C(i,j,\ell)=0$ for any $i>0,j\geq 0$,

\vspace{2mm}

$\bullet$ if $y=0$ then $C(i,j,\ell)=0$ for any $i\geq 0,j>0$.

\vspace{2mm}

Therefore in $M_{R,\ell}$ there are rows identically zero, precisely those rows corresponding to points in $R\setminus R_0$.

\vspace{2mm}

Now notice that the set $R_0$ satisfies the monomial condition.
We consider the system, for $\ell\geq r_0$,
\[M_{R_0,\ell}\big({\bf a},1\big)^t=\bar 0\]
in $\C[x,y,z]$. Since the polynomial $\prod_{(i,j)\in R_0}\big(s-A_{ij}\big)^{n_{ij}}$ is invariant under the specialization $x=0$, it is enough to apply Theorem \ref{TheSt}.

\vspace{2mm}

Finally det$\big(M_{R_0}\big)=\prod_{j=0}^{r_0-1}C(0,j,j)$ and $C(0,j,j)=y^j$.
\hfill$\Box$

\vspace{2mm}

\begin{Rem}
It is straightforward to check that the formula in Proposition \ref{Result2} holds in these cases.
\end{Rem}

\vspace{5mm}

The cases (i) and (ii) are much more involved, and are rather different. We study them separately. From now on, we are working in the ring
\[\C[x,y,z]/(ax+by)\ \mbox{ with }ab\neq 0\mbox{ and }a>0\]
unless otherwise stated.

\subsection{Systems in $\C[x,y,z]/(ax+by)$ with $b>0$}
Let us study the system $M_{R,\ell}\big({\bf a},1\big)^t=\bar 0$ in the ring $\C[x,y,z]/(ax+by)$.

First we study the specializations of the polynomials $C(i,j,\ell)$ in $\C[x,y,z]/(ax+by)$. The explicit description of $C(i,j,\ell)$ under the condition $ax+by=0$ is of no help in this case, but instead we find useful linear relations among them.

\begin{Lem}
Let $(a,b)\in\Z_{>0}^2$  and $(k_1,k_2)\in\Z_{\geq 0}^2$. In the ring $\C[x,y,z]/(ax+by)$ we have the following identities, for any $\ell\geq 0$:

\begin{equation}
%%C(k_1,k_2,l)=\sum_{i=0}^a\sum_{j=0}^b\frac{a!}{(a-i)!}\frac{b!}{(b-j)!}\binom{k_1+i}{k_1}\binom{k_2+j}{k_2}C(k_1+i,k_2+j,\ell)
\sum_{(i,j)\in S_{a,b}\setminus\{(0,0)\}}\frac{a!}{(a-i)!}\frac{b!}{(b-j)!}\binom{k_1+i}{k_1}\binom{k_2+j}{k_2}C(k_1+i,k_2+j,\ell)=0
\label{Puti2}
\end{equation}
\label{LemGP0}
\end{Lem}

{\em Proof.}  Recall that we abuse of notation by denoting $C(i,j,\ell)$ the polynomial under the condition $ax+by=0$. The proof is by induction on $k_1+k_2$. For $k_1+k_2=0$ we have that $(k_1,k_2)=(0,0)$ and, by Lemma \ref{combinatorics}, it follows that
$$\sum_{(i,j)\in S_{a,b}}\frac{a!}{(a-i)!}\frac{b!}{(b-j)!}C(i,j,\ell)=A_{ab}^\ell.$$
Then, since $A_{ab}=ax+by+z=z=A_{00}$, the first step of induction follows.

\vspace{2mm}

Suppose that the claim is true for any $(k_1,k_2)$ with $k_1+k_2\leq n$. Let us prove it for $(k_1,k_2)$ with $k_1+k_2=n+1$. We denote equation (\ref{Puti2}) as $R(k_1,k_2)=0$. We are going to prove that
$$\sum_{(i,j)\in S_{k_1,k_2}}\binom{k_1}{i}\binom{k_2}{j}i!j!R(i,j)=0,$$
and since for $(i,j)\in S_{k_1,k_2}^\bullet$ we have $R(i,j)=0$ by induction hypothesis, we deduce that $R(k_1,k_2)=0$, as we wanted to prove.

\vspace{2mm}

We have that $\sum_{(i,j)\in S_{k_1,k_2}}\binom{k_1}{i}\binom{k_2}{j}i!j!R(i,j)=$
$$=\sum_{(i,j)\in S_{k_1,k_2}}\binom{k_1}{i}\binom{k_2}{j}i!j!\left\{\sum_{(\lambda,\mu)\in S_{a,b}}\binom{a}{\l}\binom{b}{\mu}\binom{i+\l}{i}\binom{j+\mu}{j}\l!\mu!C(i+\l,j+\mu,\ell)-C(i,j,\ell)\right\}.$$
And by Lemma \ref{combinatorics} (i) this is equal to
$$\sum_{(i,j)\in S_{k_1,k_2}}\sum_{(\l,\mu)\in S_{a,b}}\binom{a}{\l}\binom{b}{\mu}\binom{k_1}{i}\binom{k_2}{j}\binom{i+\l}{i}\binom{j+\mu}{j}\l!\mu!i!j!C(i+\l,j+\mu,\ell)-A_{k_1k_2}^\ell.$$
Under the condition $ax+by=0$ we have $A_{k_1k_2}=A_{k_1+a,k_2+b}$, and hence we have to prove that
$$\sum_{(i,j)\in S_{k_1,k_2}}\sum_{(\l,\mu)\in S_{a,b}}\binom{a}{\l}\binom{b}{\mu}\binom{k_1}{i}\binom{k_2}{j}\binom{i+\l}{i}\binom{j+\mu}{j}\l!\mu!i!j!C(i+\l,j+\mu,\ell)=A_{k_1+a,k_2+b}^\ell.$$
The left hand side of the equation above is equal to
$$\sum_{(r,s)\in S_{k_1+a,k_2+b}}\left\{\sum_{i+\l=r,0\leq i\leq k_1,0\leq\l\leq a}\sum_{j+\mu=s,0\leq j\leq k_2,0\leq \mu\leq b}\binom{a}{\l}\binom{b}{\mu}\binom{k_1}{i}\binom{k_2}{j}\binom{r}{i}\binom{s}{j}\l!\mu!i!j!\right\}C(r,s,\ell)=$$
$$=\sum_{(r,s)\in S_{k_1+a,k_2+b}}\left(\sum_{r=i+\l,0\leq i\leq k_1,0\leq\l\leq a}\binom{a}{\l}\binom{r}{i}\binom{k_1}{i}\l!i!\right)\left(\sum_{s=j+\mu,0\leq j\leq k_2,0\leq \mu\leq b}\binom{b}{\mu}\binom{s}{j}\binom{k_2}{j}\mu!j!\right)C(r,s,\ell)=$$
$$=\sum_{(r,s)\in S_{k_1+a,k_2+b}}\left(r!\sum_{r=i+\l,0\leq i\leq k_1,0\leq\l\leq a}\binom{a}{\l}\binom{k_1}{i}\right)\left(s!\sum_{s=j+\mu,0\leq j\leq k_2,0\leq\mu\leq b}\binom{b}{\mu}\binom{k_2}{j}\right)C(r,s,\ell).$$
By Vandermonde Convolution (\ref{VDMD}), we have that this is equal to
$$\sum_{(r,s)\in S_{k_1+a,k_2+b}}r!\binom{k_1+a}{r}s!\binom{k_2+b}{s}C(r,s,\ell)$$
which, by Lemma \ref{combinatorics} (i), is equal to $A_{k_1+a,k_2+b}^\ell$, and this finishes the proof.
\hfill$\Box$

\vspace{3mm}

\begin{Defi}
Given a set of points $R$ satisfying the monomial condition, and $(a,b)\in\Z_{>0}^2$, we define the subset of $R$
\[R_{a,b}=\{(i,j)\in R\ |\ (i,j)-(a,b)\notin R\}.\]
\label{eqrel0}
\end{Defi}

\begin{Rem}
We have that $R_{a,b}=R$ unless there are $(i,j),(i'j')\in R$ with
\[(i',j')-(i,j)=(a,b).\]
Note also that $R_{a,b}=R$ is equivalent to $b_R(s)$ non-reduced in $\C[x,y,z]/(ax+by)$.
\label{RemJunio}
\end{Rem}

\begin{Pro}
For any $\ell\in\Z_{\geq 0}$ the system $M_{R,\ell}\big({\bf a},1\big)^t=\bar 0$ is equivalent to the system $M_{R_{a,b},\ell}\big({\bf a},1\big)^t=\bar 0$ in $\C[x,y,z]/(ax+by)$.
\end{Pro}

{\em Proof.}
Note that Lemma \ref{LemGP0} provides relations among the rows of the matrix $M_{R,\ell}$ in $\C[x,y,z]/(ax+by)$.
We have to prove that for any $(k_1,k_2)\in R$ such that $(k_1,k_2)-(a,b)\in R$, we can drop the row corresponding to $(k_1,k_2)$ in the system $M_{R,\ell}\big({\bf a},1\big)^t=\bar 0$. This follows by Lemma \ref{LemGP0}, using that, $S_{k_1,k_2}\subseteq R$ for any $(k_1,k_2)\in R$, since $R$ has the monomial condition.
\hfill$\Box$

\vspace{3mm}

\begin{Lem}
If $R$ satisfies the monomial condition, then, for any $(a,b)\in\Z_{>0}^2$, the set
$R_{a,b}$ also satisfies the monomial condition.
\end{Lem}

{\em Proof.} It is enough to prove that if $(k_1,k_2)\in R_{a,b}$ then $(k_1-1,k_2),(k_1,k_2-1)\in R_{a,b}$.

If $(k_1-1,k_2)\in\Z_{\geq 0}^2$ but $(k_1-1,k_2)\notin R_{a,b}$ then
\[(k_1-1,k_2)-(a,b)\in R.\]
But then, since $R$ has the monomial condition, we deduce that $(k_1,k_2)-(a,b)\in R$, which contradicts that $(k_1,k_2)\in R_{a,b}$.

Analogously, if $(k_1,k_2-1)\notin R_{a,b}$, then $(k_1,k_2-1)-(a,b)\in R$, but then, since $S_{k_1,k_2}\subset R$ (because $R$ has the monomial condition), we deduce that  $(k_1,k_2)-(a,b)\in R$.
\hfill$\Box$

\vspace{3mm}

\begin{Pro}
Given a set of $r$ points in $\Z_{\geq 0}^2$ of the form $R_{a,b}$ where $R$ satisfies the monomial condition, and $(a,b)\in\Z_{>0}^2$, consider the system
\[M_{R_{a,b},\ell}\big({\bf a},1\big)^t=\bar 0\]
with $\ell\geq r$. In the ring $\C[x,y,z]/(ax+by)$, any ${\bf a}=(a_0,\ldots,a_{\ell-1})$ such that $b_{\bf a}(s)=\prod_{(i,j)\in R_{a,b}}(s-A_{ij})^{n_{ij}}$ with $n_{ij}>0$ and $\sum n_{ij}=\ell$, is a solution to the system.
\label{TheP0}
\end{Pro}

{\em Proof.} It is a consequence of Theorem \ref{TheSt}, since $R_{a,b}$ has the monomial condition.\hfill$\Box$

\vspace{3mm}

\begin{Rem}
Since $R_{a,b}$ satisfies the monomial condition, we can apply Proposition \ref{determinant} to compute det$(M_{R_{a,b}})$. By definition of $R_{a,b}$ and Remark \ref{RemJunio}, we deduce that det$(M_{R_{a,b}})$ is non-zero in the ring $\C[x,y,z]/(ax+by)$.

\end{Rem}

\vspace{3mm}

\subsection{Systems in $\C[x,y,z]/(ax-by)$ with $b>0$}
For convenience we write the condition $ax+by=0$ with $b<0$ as $ax-by=0$ with $b>0$. As in the previous case there are relations among the polynomials $C(i,j,\ell)$ in the ring $\C[x,y,z]/(ax-by)$.

\begin{Lem}
Let $(a,b)\in\Z_{>0}^2$  and $(k_1,k_2)\in\Z_{\geq 0}^2$ such that $k_2-b\geq 0$. In the ring $\C[x,y,z]/(ax-by)$ we have the following identities for any $\ell\geq 0$:
\begin{equation}
\sum_{j=1}^a\frac{a!}{(a-j)!}\binom{k_1+j}{k_1}C(k_1+j,k_2-b,\ell)=\sum_{j=1}^{b}\frac{b!}{(b-j)!}\binom{k_2-b+j}{k_2-b}C(k_1,k_2-b+j,\ell)
\label{GranPuta}
\end{equation}
\label{LemGP}
\end{Lem}

{\em Proof.} We abuse of notation by denoting $C(i,j,\ell)$ the polynomial under the condition $ax=by$.
 We will prove the statement by induction on $k_1+k_2\geq b$.

The first step of induction corresponds to the point $(k_1,k_2)=(0,b)$. Under the condition  $ax-by=0$ we obviously have $A_{a0}^\ell=A_{0b}^\ell$ for any $\ell$. Then by Lemma \ref{combinatorics} (i) we deduce that, for any $\ell\geq 0$
\[\sum_{i=0}^a\binom{a}{i}i!{C}(i,0,\ell)=\sum_{j=0}^{b}\binom{b}{j}j!{C}(0,j,\ell).\]
Notice that $C(0,0,\ell)$ cancels on the previous equality and then we get the claim.

\vspace{2mm}

Suppose now that the claim is true for any $(k_1,k_2)$ with $k_2\geq b$ and $k_1+k_2\leq b+r$ and we will prove it for $(k_1,k_2)$ with $k_2\geq b$ and $k_1+k_2=b+r+1$.
For any $(i,j)$ with $j\geq b$ and $i+j\leq b+r$ we have the relation in the statement, by induction hypothesis. Let us denote by $R(i,j)=0$ such relation.
Then, for any $(i,j)\in S^\bullet_{k_1,k_2}\setminus S_{k_1,b-1}$, we can apply induction hypothesis to get
$$R(i,j)=0.$$
Hence, if we prove that
\begin{equation}
\sum_{(i,j)\in S_{k_1,k_2}\setminus S_{k_1,b-1}}\binom{k_1}{i}\binom{k_2-b}{j-b}i!(j-b)!R(i,j)=0
\label{Srelation}
\end{equation}
we deduce that $R(k_1,k_2)=0$, as we wanted to prove.

\vspace{2mm}

Let us prove equation (\ref{Srelation}). First we deal with the part:
\[\sum_{(i,j)\in S_{k_1,k_2}\setminus S_{k_1,b-1}}\binom{k_1}{i}\binom{k_2-b}{j-b}i!(j-b)!\left\{\sum_{r=1}^a\binom{a}{r}r!\binom{i+r}{i}{C}(i+r,j-b,\ell)\right\}.\]
The sum above is equal to
\[\sum_{(i,j)\in S_{k_1,k_2}\setminus S_{k_1,b-1}}i!(j-b)!\sum_{r=0}^a\binom{a}{r}\binom{k_1}{i}\binom{k_2-b}{j-b}\binom{i+r}{i}r!{C}(i+r,j-b,\ell)-\]
\[-\sum_{(i,j)\in S_{k_1,k_2}\setminus S_{k_1,b-1}}i!(j-b)!\binom{a}{0}\binom{k_1}{i}\binom{k_2-b}{j-b}{C}(i,j-b,\ell)=\]
\[\sum_{(i,j)\in S_{k_1,k_2-b}}i!j!\sum_{r=0}^a\binom{a}{r}\binom{k_1}{i}\binom{k_2-b}{j}\binom{i+r}{i}r!{C}(i+r,j,\ell)-\sum_{(i,j)\in S_{k_1,k_2-b}}i!j!\binom{k_1}{i}\binom{k_2-b}{j}{C}(i,j,\ell)=\]
\[=\sum_{(\lambda,\mu)\in S_{k_1+a,k_2-b}}\left\{\sum_{0\leq i\leq k_1,0\leq r\leq a,i+r=\lambda}i!\mu!\binom{a}{r}\binom{k_1}{i}\binom{k_2-b}{\mu}\binom{i+r}{i}r!\right\}{C}(\lambda,\mu,\ell)-\]
\[-\sum_{(i,j)\in S_{k_1,k_2-b}}i!j!\binom{k_1}{i}\binom{k_2-b}{j}{C}(i,j,\ell)=\]
\[=\sum_{(\lambda,\mu)\in S_{k_1+a,k_2-b}}\binom{k_2-b}{\mu}\mu!\left\{\sum_{0\leq i\leq k_1,0\leq r\leq a,i+r=\lambda}i!\binom{a}{r}\binom{k_1}{i}\binom{i+r}{i}r!\right\}{C}(\lambda,\mu,\ell)-\]
\[-\sum_{(i,j)\in S_{k_1,k_2-b}}i!j!\binom{k_1}{i}\binom{k_2-b}{j}{C}(i,j,\ell)=\]
\[=\sum_{(\lambda,\mu)\in S_{k_1+a,k_2-b}}\binom{k_2-b}{\mu}\mu!\lambda!\left\{\sum_{0\leq i\leq k_1,0\leq r\leq a,i+r=\lambda}\binom{a}{r}\binom{k_1}{i}\right\}{C}(\lambda,\mu,\ell)-\]
\[-\sum_{(i,j)\in S_{k_1,k_2-b}}i!j!\binom{k_1}{i}\binom{k_2-b}{j}{C}(i,j,\ell).\]
And this is equal to
\[\sum_{(\lambda,\mu)\in S_{k_1+a,k_2-b}}\binom{k_1+a}{\lambda}\binom{k_2-b}{\mu}\mu!\lambda!{C}(\lambda,\mu,\ell)-\sum_{(i,j)\in S_{k_1,k_2-b}}i!j!\binom{k_1}{i}\binom{k_2-b}{j}{C}(i,j,\ell),\]
by Vandermonde Convolution (\ref{VDMD}).

\vspace{3mm}

Analogously the second part of the relation we want to prove is
\[\sum_{(i,j)\in S_{k_1,k_2}\setminus S_{k_1,b-1}}i!(j-b)!\binom{k_1}{i}\binom{k_2-b}{j-b}\sum_{r=1}^{b}\binom{b}{r}r!\binom{j-b+r}{j-b}{C}(i,j-b+r,\ell)=\]
\[=\sum_{(\lambda,\mu)\in S_{k_1,k_2}}\binom{k_1+a}{\lambda}\binom{k_2-b}{\mu}\lambda!\mu!{C}(\lambda,\mu,\ell)-\sum_{(i,j)\in S_{k_1,k_2-b}}i!j!\binom{k_1}{i}\binom{k_1-b}{j}{C}(i,j,\ell).\]

Hence equation (\ref{Srelation}) is equivalent to
\[\sum_{(\lambda,\mu)\in S_{k_1+a,k_2-b}}\binom{k_1+a}{\lambda}\binom{k_2-b}{\mu}\lambda!\mu!{C}(\lambda,\mu,\ell)=\sum_{(\lambda,\mu)\in S_{k_1,k_2}}\binom{k_1+a}{\lambda}\binom{k_2-b}{\mu}\lambda!\mu!{C}(\lambda,\mu,\ell).\]
And we are done, because this equality holds by Lemma \ref{combinatorics} (i), since under the condition $ax-by=0$ we have $A_{k_1,k_2}=A_{k_1+a,k_2-b}$.
\hfill$\Box$

\vspace{5mm}

Let us define the following equivalence relation on $\Z^2$:
\begin{equation}
(i,j)\sim_{a,-b}(i',j')\mbox{ if and only if }(i,j)=(i',j')+\lambda(a,-b)
\label{defRE}
\end{equation}
for $\lambda\in\Q$, and consider the quotient space
\[R/\sim_{a,-b}.\]
Since in the ring $\C[x,y,z]/(ax-by)$ we have
\[\big(b_R(s)\big)_{red}=b_{R/\sim_{a,-b}}(s),\]
it seems plausible that $R/\sim_{a,-b}$ is the set we are looking for.

Note that the set $R/\sim_{a,-b}$ can be identified with the set
\[\{(i,j)\in R\ |\ (i,j)+(a,-b)\notin R\}.\]
Indeed, this corresponds with choosing, for each equivalence class, the representant of the smallest height, and it is exactly what we did in the previous case (without mention any equivalence relation). The problem is that in general it does not satisfy the monomial condition (see Example \ref{ejemplito}). Is there any other choice of representant such that $R/\sim_{a,-b}$ inherits form $R$ the monomial condition? Unfortunately not, as we see in the following example.

\vspace{2mm}

\begin{Exam}
Consider the set $R=\{(0,0),(1,0),(2,0),(3,0),(0,1),(0,2),(0,3),(0,4)\}$ and the equivalence relation $\sim_{2,-3}$. The only choices of $R/\sim_{2,-3}$ are
\[\{(0,0),(1,0),(2,0),(3,0),(0,1),(0,2),(0,4)\}\]
\[\{(0,0),(1,0),(3,0),(0,1),(0,2),(0,3),(0,4)\}\]
and none of them satisfies the monomial condition.
\label{ejemplito}
\end{Exam}

\vspace{2mm}

Let us study then the set $R/\sim_{a,-b}$ even though we know that it does not satisfy the monomial condition. Observe that in this case there are many candidates for $R_{a,b}$: in fact every choice of representant in $R/\sim_{a,-b}$ will work.

Note that though the polynomial $b_{R/\sim_{a,-b}}(s)$ is well defined in the ring $\C[x,y,z]/(ax-by)$, the system $M_{R/\sim_{a,-b},\ell}\big({\bf a},1\big)^t=\bar 0$ does depend on the choice of representant we make on each equivalence class. We will use Lemma \ref{LemGP} to define the system $M_{R/\sim_{a,-b},\ell}\big({\bf a},1\big)^t=\bar 0$ properly in  $\C[x,y,z]/(ax-by)$.

\vspace{3mm}

\begin{Pro}
Let $R\subseteq\Z^2$ be a set satisfying the monomial condition. For any $(a,b)\in\Z_{>0}^2$ and $\ell\in\Z_{\geq 0}$ the system
 $M_{R,\ell}({\bf a},1)^t=\overline 0$ in the ring $\C[x,y,z]/(ax-by)$ is equivalent to
\[M_{R/\sim_{a,-b},\ell}({\bf a},1)^t=\overline 0\]
for any choice of representant in $R/\sim_{a,-b}$. Recall that ${\bf a}=(a_0,\ldots,a_{\ell-1})$.
\label{CorST}
\end{Pro}

{\em Proof.}
Notice that whenever $(k_1,k_2)$ and $(k_1+a,k_2-b)$ are points in $R$, we have that
\[(k_1+j,k_2-b)\in R\mbox{ for }j=1,\ldots,a,\]
\[(k_1,k_2-b+j)\in R\mbox{ for }j=1,\ldots,b,\]
since $R$ has the monomial condition.

\vspace{2mm}

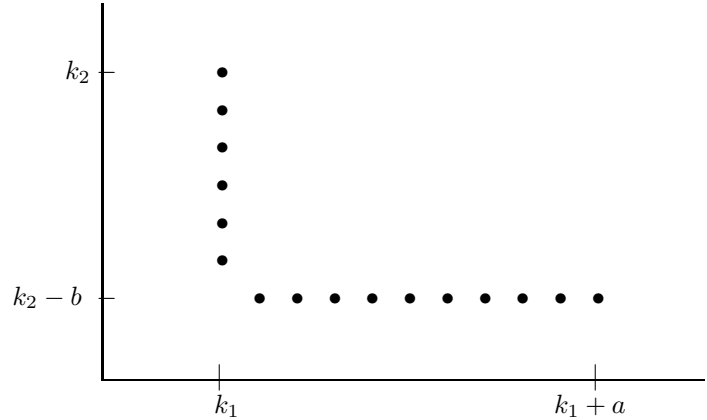
\begin{figure}[h]
\unitlength=1mm
\centering
\begin{center}
\begin{picture}(70,48)
\linethickness{0.15mm}
\put(0,0){\line(1,0){80}}
\put(0,0){\line(0,1){50}}
\put(-0.8,10){$-$}\put(-12,10){$k_2-b$}
\put(-0.8,40){$-$}\put(-5,40){$k_2$}
\put(15,-0.5){$|$}\put(14.8,-4.4){$k_1$}
\put(65,-0.5){$|$}\put(60,-4.4){$k_1+a$}
\put(15,40){$\bullet$}\put(60,10){$\bullet$}\put(65,10){$\bullet$}
\put(15,35){$\bullet$}\put(55,10){$\bullet$}
\put(15,30){$\bullet$}\put(50,10){$\bullet$}
\put(15,25){$\bullet$}\put(45,10){$\bullet$}
\put(15,20){$\bullet$}\put(40,10){$\bullet$}
\put(15,15){$\bullet$}\put(35,10){$\bullet$}
\put(30,10){$\bullet$}\put(25,10){$\bullet$}\put(20,10){$\bullet$}
\end{picture}
\end{center}
\vspace{3mm}
\caption{If $ax-by=0$ with $b>0$, then, for any $\ell\in\Z_{\geq 0}$ we have a relation among the $C(i,j,\ell)$ corresponding to the points represented in the picture.}
\label{FNoro}
\end{figure}

 By Lemma \ref{LemGP} there is a linear combination among the rows of the matrix $M_{R,\ell}$ corresponding to the points
\begin{equation}
\begin{array}{l}
(k_1+1,k_2-b),\ldots,(k_1+a,k_2-b)\\
\\
(k_1,k_2-b+1),\ldots,(k_1,k_2)\\
\end{array}
\label{PuntosP}
\end{equation}
(see Figure \ref{FNoro}).
If $(k_1,k_2),(k_1+a,k_2-b)\in R$, by the monomial condition we deduce that all the points in (\ref{PuntosP}) are in the set $R$. By Lemma \ref{LemGP} we can drop a row corresponding to one of the points in (\ref{PuntosP}), obtaining in this way an equivalent system of equations.

Now we repeat the argument as long as there exist points of the form $(i,j)$ and $(i,j)+\lambda(a,-b)$ in the remaining set of points. Note that now it may be that not all the points in (\ref{PuntosP}) are in the set, but if one is missing, the reason is because we have taken it away in a previous step, and then we can write the corresponding row as a linear combination of other rows.

At any step we decide which row to drop. This corresponds with deciding which element represents a class in $R/\sim_{a,-b}$. Dropping every element on a class except one, we prove that the system $M_{R,\ell}\big({\bf a},1\big)^t=\bar 0$ is equivalent to $M_{R/\sim_{a,-b},\ell}\big({\bf a},1\big)^t=\bar 0$ for any $\ell$.

Note that this argument also works with the points $(k_1,k_2)$ and $(k_1+na,k_2-nb)$, since in Lemma \ref{LemGP} we do not ask $a$ and $b$ to be coprime.
\hfill$\Box$

\vspace{3mm}

\begin{Rem}
Note that this is not true in the ring $\C[x,y,z]/(ax+by)$ with $b>0$. If we consider the equivalence relation analogously as in (\ref{defRE}), not every choice of representant yields to equivalent systems. In general only the set $R_{a,b}$ defined in Definition \ref{eqrel0} works, which corresponds to choosing the representant with the smallest height.
\end{Rem}

\vspace{2mm}

Now we prove an analogous result to Theorem \ref{TheSt} in the ring $\C[x,y,z]/(ax-by)$. The proof is not as direct as the proof of Proposition \ref{TheP0}, because  $R/\sim_{a,-b}$ does not have the monomial condition.

\begin{The}
Consider a set of $r$ points in $\Z_{\geq 0}^2$ of the form $R/\sim_{a,-b}$ where $R$ satisfies the monomial condition, and $(a,b)\in\Z_{>0}^2$.
 For any $\ell\geq r$, we have that any ${\bf a}=(a_0,\ldots,a_{\ell-1})$ such that $b_{\bf a}(s)=\prod_{(i,j)\in R/\sim_{a,-b}}(s-A_{ij})^{n_{ij}}$ with $n_{ij}>0$ and $\sum n_{ij}=\ell$, is a solution to the system.
\[M_{R/\sim_{a,-b},\ell}\big({\bf a},1\big)^t=\bar 0,\]
in the ring $\C[x,y,z]/(ax-by)$.
\label{TheP}
\end{The}

{\em Proof.} It is enough to prove it for a particular choice in $R/\sim_{a,-b}$, since, by Proposition \ref{CorST}, all these systems are equivalent. For any class in $R/\sim_{a,-b}$ we choose the representant with the smallest height.

The proof is analogous to the proof of Theorem \ref{TheSt}, by double induction  on $r>0$ and $\ell\geq r$.

\vspace{2mm}

For $r=1$ the set $R/\sim_{a,-b}$ is the set $\{(0,0)\}$, since $(0,0)\not\sim_{a,-b}(i,j)\in\Z_{\geq 0}^2$, and the proof is exactly as in the proof of Theorem \ref{TheSt}.

For $r>1$, suppose first that $\ell=r$. We only need to prove that we can always find $(k_1,k_2)\in R/\sim_{a,-b}$ such that
\[R/\sim_{a,-b}\setminus\{(k_1,k_2)\}=R'/\sim_{a,-b}\]
where $R'$ has the monomial condition, so that we can apply induction exactly as in the proof of Theorem \ref{TheSt}. Indeed, notice that on each equivalence class in $R/\sim_{a,-b}$ the quantity $bi+aj$ is constant. Consider the equivalence class which maximizes $bi+aj$. Then the representant $(k_1,k_2)$ with the smallest height belongs to our choice of $R/\sim_{a,-b}$, the set
\[R'=R\setminus\{(i,j)\in R\ |\ (i,j)\sim_{a,-b}(k_1,k_2)\}\]
has the monomial condition and $R/\sim_{a,-b}\setminus\{(k_1,k_2)\}=R'/\sim_{a,-b}$. To follow the proof of Theorem \ref{TheSt} note that in general $S_{k_1,k_2}^\bullet\nsubseteq R'/\sim_{a,-b}$, but we can apply Lemma \ref{LemGP} and deduce the equality (\ref{july}) in this case.

The proof for $\ell>r$ goes exactly as in Theorem \ref{TheSt}.\hfill$\Box$

\vspace{3mm}

 We finish this section by studying the determinant of the matrix $M_{R/\sim_{a,-b}}$. In this case we can not apply Proposition \ref{determinant}, since in general $R/\sim_{a,-b}$ does not have the monomial condition. However we have the following analogous result.

\begin{Pro}
 Let $H$ be an ordered set of points in $\Z^2$ consisting of a particular choice of representant in the quotient $R/\sim_{a,-b}$, for certain $(a,b)\in\Z^2_{0}$, where $R$ is a set of points satisfying the monomial condition. Then, in the ring $\C[x,y,z]/(ax-by)$ we have
\[det(M_H)=\prod_{(k_1,k_2)\in H}\frac{1}{k_1!k_2!}\prod_{(i,j)<(i',j')\in H}(A_{ij}-A_{i'j'})\]
where the order $<$ is given by the order of the points in the rows of the matrix.
\end{Pro}

{\em Proof.}
Exactly as in the proof of Proposition \ref{determinant} we have
\[\mbox{det}(M_H)=\prod_{(k_1,k_2)\in H}\frac{1}{k_1!k_2!}\mbox{ det}\left(\left\langle\frac{t^\ell}{\ell!}\right\rangle(e^{xt}-1)^{k_1}(e^{yt}-1)^{k_2}e^{zt}\right)_{(k_1,k_2)\in H,\ 0\leq \ell<r}\]
Let $r$ be the cardinal of the set $H$. If the points in $H$ are ordered as $\{(k_1^{(1)},k_2^{(1)}),\ldots,(k_1^{(r)},k_2^{(r)})\}$, let us define the set
\[\mathcal S_H=S_{k_1^{(1)},k_2^{(1)}}\times S_{k_1^{(2)},k_2^{(2)}}\times\cdots\times S_{k_1^{(r)},k_2^{(r)}}.\]
Hence, since for any $\ell$
\[\left\langle\frac{t^\ell}{\ell!}\right\rangle(e^{xt}-1)^{k_1}(e^{yt}-1)^{k_2}e^{zt}=\left\langle\frac{t^\ell}{\ell!}\right\rangle\sum_{(i,j)\in S_{k_1,k_2}}\binom{k_1}{i}\binom{k_2}{j}(-1)^{k_1+k_2-i-j}e^{(ix+jy+z)t},\]
we have that
\[\begin{array}{ll}
D_H & :=\mbox{det}\left(\left\langle\frac{t^\ell}{\ell!}\right\rangle(e^{xt}-1)^{k_1}(e^{yt}-1)^{k_2}e^{zt}\right)_{(k_1,k_2)\in H,\ 0\leq \ell<r}\\
\\
    & =\sum_{s\in\mathcal S_H}\prod_{q=1}^r\binom{k_1^{(q)}}{i^{(q)}}\binom{k_2^{(q)}}{j^{(q)}}(-1)^{k_1^{(q)}+k_2^{(q)}-i^{(q)}-j^{(q)}}\mbox{det}(\Delta_s)\\
\end{array}\]
where the tuple $s=(i^{(1)},j^{(1)},i^{(2)},j^{(2)},\ldots,i^{(r)},j^{(r)})\in\mathcal S_H$ will be identified when necessary with the ordered set $$\{(i^{(1)},j^{(1)}),\ldots,(i^{(r)},j^{(r)})\},$$ and $\Delta_s$ is the Vandermonde matrix
\[\Delta_s=\left( \big(i^{(q)}x+j^{(q)}y+z\big)^\ell\right)_{1\leq q\leq r,\ 0\leq \ell<r}\]
Then
\[\mbox{det}(\Delta_s)=\prod_{1\leq p<q\leq r}\big(A_{i^{(p)}j^{(p)}}-A_{i^{(q)}j^{(q)}}\big)\]
We clearly have
\[H\in\mathcal S_H\]
and
\[\mathcal S_H\subseteq\bigcup_{(k_1,k_2)\in H}S_{k_1,k_2}\subseteq R\]
because $R$ satisfies the monomial property.

\vspace{2mm}

We claim that for all $s\in\mathcal S_H\setminus\{H\}$
\[\mbox{det}(\Delta_s)=0\ \mbox{ in }\ \C[x,y,z]/(ax-by)\]
and this proves the statement. Let us prove the claim.
If there are repeated points in $s$, we clearly have
\[\mbox{det}(\Delta_s)=0\ \mbox{ in }\C[x,y,z]\]
Otherwise, $s$ describes $r$ different points in $$\bigcup_{(k_1,k_2)\in H}S_{k_1,k_2}\subseteq R.$$ If $s/\sim_{a,-b}$ has cardinal smaller than $r$, then there are at least two points related by $\sim_{a,-b}$, or equivalently $A_{i^{(p)}j^{(p)}}\equiv A_{i^{(q)}j^{(q)}}\mbox{ mod }ax-by$ for certain $p$ and $q$, and the determinant vanishes in the ring $\C[x,y,z]/(ax-by)$. The remaining case is when $s/\sim_{a,-b}$ consists of $r$ equivalence classes. Then $s$ must be another choice of representant in $R/\sim_{a,-b}$, because $s\in R$. But any choice of $R/\sim_{a,-b}$ different from $H$ can not be a point in $\mathcal S_H$. Indeed, reordering the elements if necessary we can suppose that the points are ordered in increasing weight $bk_1+ak_2$ (note that the equivalence classes of $R/\sim_{a,-b}$ are in correspondence with the numbers $bk_1+ak_2$ for $(k_1,k_2)\in H$). Let $q\in\{1,\ldots,r\}$ be maximum with the property
\[(i^{(q)},j^{(q)})\in S_{k_1^{(q)},k_2^{(q)}}^\bullet\]
Such a maximum exists because $s\in\mathcal S_H\setminus \{H\}$. Then
\[(i^{(q)},j^{(q)})\in\bigcup_{p=1}^{q-1}S_{k_1^{(i)},k_2^{(i)}}\]
but this is impossible, since any $(i,j)\in S_{k_1,k_2}$ has weight $bi+aj\leq bk_1+ak_2$ . Therefore there are no points of weight $bi^{(q)}+aj^{(q)}$  in the set $\cup_{p=1}^{q-1}S_{k_1^{(i)}k_2^{(i)}}$.
\hfill$\Box$

\vspace{3mm}

\begin{Rem}
Notice that, by definition of the relation $\sim_{a,-b}$, for any choice $H$ of representant in $R/\sim_{a,-b}$, there are no points $(i,j),(i',j')\in H$
such that $A_{ij}-A_{i'j'}=(i-i')x+(j-j')y\equiv 0\mbox{ mod }(ax-by)$, and hence det$(M_H)\neq 0$ in $\C[x,y,z]/(ax-by)$.
\end{Rem}

\vspace{3mm}

\section{Some combinatorial applications}
\label{Appl}
In this section we highlight some relations on Stirling numbers of second kind that were developed in previous sections. We have not found any reference for these relations.

\vspace{2mm}

\begin{enumerate}
\item[(i)]  Notice that, with our notations, the formula in (\ref{eqPalma}) can be written as
\[\sum_{k=0}^n\binom{n}{k}(-1)^kA_{k0}^m=(-1)^nn!C(n,0,m),\]
which looks like a particular case of the equality in Lemma \ref{combinatorics} (ii). Indeed, dropping our notation, this equality can be written as
\[\sum_{(i,j)\in S_{k_1,k_2}}\frac{(-1)^{k_1+k_2-i-j}}{i!j!(k_1-i)!(k_2-j)!}\big(ix+jy+z\big)^\ell=\sum_{i_1\geq k_1,i_2\geq k_2,i_1+i_2\leq\ell}\binom{\ell}{i_1}\binom{\ell-i_1}{i_2}S(i_1,k_1)S(i_2,k_2)x^{i_1}y^{i_2}z^{\ell-i_1-i_2}\]
which is a generalization of equation (\ref{eqPalma}).

\vspace{5mm}

\item[(ii)]  We have given in Lemma \ref{LemGP0} and Lemma \ref{LemGP} explicit linear relations among the rows of the matrix $M_R$ under some particular specializations. As a consequence we will recover interesting relations among the Stirling numbers of second kind. Let us see some examples.

\begin{Exam}
If $x=y$ and $z=0$, we have that, for any $(k_1,k_2)\in\Z^2_{\geq 0}$ with $k_2>0$,
\[C(k_1,k_2,\ell)=x^\ell\sum_{i_1=k_1}^{\ell-k_2}\binom{\ell}{i_1}S(i_1,k_1)S(\ell-i_1,k_2)\]
Applying Lemma \ref{LemGP} repeatedly we have, for any $\ell$,

\[C(k_1,k_2,\ell)=\frac{k_1+1}{k_2}C(k_1+1,k_2-1,\ell)=\frac{(k_1+1)(k_1+2)}{k_2(k_2-1)}C(k_1+2,k_2-2,\ell)=\cdots=\binom{k_1+k_2}{k_1}C(k_1+k_2,0,\ell)\]
which yields to the well-known convolution relation
\[\binom{k_1+k_2}{k_1}S(\ell,k_1+k_2)=\sum_{i_1=k_1}^{\ell-k_2}\binom{\ell}{i_1}S(\ell-i_1,k_2)S(i_1,k_1)\]
\end{Exam}

\begin{Exam}
If $b=1$, $y=ax$ and $z=0$, we deduce, by Lemma \ref{LemGP} applied to the point $(k_1,1)$, that
\[\sum_{j=1}^a\frac{a!}{(a-j)!}\binom{k_1+j}{k_1}S(\ell,k_1+j)=\sum_{i=k_1}^{\ell-1}\binom{\ell}{i}S(i,k_1)a^{\ell-i}\]
While if we apply it to the point $(0,k_2)$ we deduce
\[k_2S(\ell,k_2)=\sum_{j=1}^a\frac{a!}{(a-j)!}\sum_{i=j}^{\ell-k_2+1}\binom{\ell}{i}\frac{1}{a^i}S(i,j)S(\ell-i,k_2-1)\]
\end{Exam}

\vspace{3mm}

\begin{Exam}
If $ax=by$ and $x=bz$ we deduce from Lemma \ref{LemGP} the relation
\[\sum_{j=1}^a\frac{a!}{(a-j)!}\sum_{i=j}^\ell\binom{\ell}{i}b^iS(i,j)=\sum_{j=1}^b\frac{b!}{(b-j)!}\sum_{i=j}^\ell\binom{\ell}{i}a^iS(i,j)\]
for any $(a,b)\in\Z^2_{>0}$.

\vspace{2mm}

And more generally, if $ax=by$ and $x=tz$ with $t\in\R\setminus\{0\}$, we have
\[\sum_{j=1}^a\frac{a!}{(a-j)!}\sum_{i=j}^\ell\binom{\ell}{i}S(i,j)t^i=\sum_{j=1}^b\frac{b!}{(b-j)!}\sum_{i=j}^\ell\binom{\ell}{i}S(i,j)\big(\frac{a}{b}t\big)^i\]
\end{Exam}

\vspace{3mm}

\item[(iii)] Both types of Stirling numbers have been generalized in many different ways (see \cite{Gen}). Many of these generalizations of the Stirling numbers of second kind are related to the polynomials $C(k_1,k_2,\ell)$, due to the relation
\[\big(e^{xt}-1\big)^{k_1}\big(e^{yt}-1\big)^{k_2}e^{zt}=k_1!k_2!\sum_{\ell=0}^\infty C(k_1,k_2,\ell)\frac{t^\ell}{\ell!}.\]
Notice that this is a direct consequence of (\ref{C2def}).

\vspace{3mm}

Let us see some examples.

\vspace{2mm}

$\bullet$    In \cite{C} the author defines the weighted Stirling numbers of second kind (see also \cite{PTY}), $S(n,w,x)$, as
\[\big(e^t-1\big)^we^{xt}=w!\sum_{n=w}^\infty S(n,w,x)\frac{t^n}{n!}\]
Notice that
\[S(n,w,z)=C(w,0,n)_{|_{x=1}}\]
and, by Lemma \ref{combinatorics} (ii) we deduce
\[n!S(n,w,x)=\sum_{i=0}^n\binom{n}{i}(-1)^{n-1}(x+i)^n\]

\vspace{3mm}

$\bullet$ In \cite{Singh} the generalized Stirling numbers of second kind are defined as
\[S^\alpha(n,k,r)=\sum_{i=0}^n\binom{n}{i}\alpha^{n-i}r^iS(i,k)\]
and this is, adapting to our notation,
\[C(k,0,n)=S^z(n,k,x).\]
By Lemma \ref{combinatorics} (ii),
\[S^z(n,k,x)=\frac{1}{k!}\sum_{i=0}^k\binom{k}{i}(-1)^{k-i}(ix+z)^n\]
Moreover, applying Lemma \ref{LemGP} (with $b=k_2=1$ and $k_1=0$) we deduce that for every $a\in\Z_{>0}$,
\[\sum_{j=1}^a\frac{a!}{(a-j)!}S^z(\ell,j,x)=(ax+z)^\ell-z^\ell\]
\end{enumerate}

\vspace{5mm}

\bibliographystyle{amsplain}
\def\cprime{$'$}
\providecommand{\bysame}{\leavevmode\hbox to3em{\hrulefill}\thinspace}
\providecommand{\MR}{\relax\ifhmode\unskip\space\fi MR }
\providecommand{\MRhref}[2]{%
  \href{http://www.ams.org/mathscinet-getitem?mr=#1}{#2}
}
\providecommand{\href}[2]{#2}

\end{document}